\documentclass[10pt,a4paper]{article}
\usepackage{amsfonts}
\usepackage{amssymb}
\usepackage{mathrsfs}
\usepackage{amsthm}
\usepackage{setspace}

\usepackage{fancyhdr}
\newtheorem{defi}{Definition}[section]
\newtheorem{theo}[defi]{Theorem}
\newtheorem{prop}[defi]{Proposition}
\newtheorem{cor}[defi]{Corollary}
\newtheorem{lemm}[defi]{Lemma}

\theoremstyle{definition}
\newtheorem*{rem}{Remark}

\newcommand{\be}{\begin{eqnarray*}}
\newcommand{\ee}{\end{eqnarray*}}
\newcommand{\beqa}{\begin{eqnarray}}
\newcommand{\eeqa}{\end{eqnarray}}
\newcommand{\ba}{\begin{array}}
\newcommand{\ea}{\end{array}}
\newcommand{\onab}{\overrightarrow{\nabla}}
\newcommand{\rP}{\mathsf{P}}

\newenvironment{lproof}{\emph{Proof of Lemma.}}{ \qed \par}

\begin{document}

\title{Projective Holonomy I: Principles and Properties}
\author{Stuart Armstrong}
\date{12 September 2006}
\maketitle

\begin{abstract}
The aim of this paper and its sequel is to introduce and classify the holonomy algebras of the projective Tractor connection. After a brief historical background, this paper presents and analyses the projective Cartan and Tractor connections, the various structures they can preserve, and their geometric interpretations. Preserved subbundles of the Tractor bundle generate foliations with Ricci-flat leaves. Contact- and Einstein-structures arise from other reductions of the Tractor holonomy, as do $U(1)$ and $Sp(1, \mathbb{H})$ bundles over a manifold of smaller dimension.
\end{abstract}

\section{Introduction}
This paper is a very small part of the ongoing effort, at least as far back as Cartan and Weyl, to attempt to put all geometries under one unifying roof -- and, just as rapidly, to cut up that roof into separate results for specific geometric structures.

The aim of this paper is to continue the project started in \cite{mecon}, that of exploring and classifying the holonomy algebras of various parabolic geometries. Papers \cite{mecon} and \cite{me!} study conformal holonomies, this one and its sequel \cite{mepro2} are interested in projective ones. Recall that a projective structure is given by a set of unparameterised geodesics, see definition \ref{defi:pro}.

Both conformal and projective geometry are members of the class of parabolic geometries, a group that includes, amongst others, almost Grassmanian, almost quaternionic, and co-dimension one CR structures. The central concepts emerged from E. Cartan's work  \cite{ECC}, \cite{ECP}, (refined with discussions and arguments with H. Weyl \cite{Weyl}), whose technique of `moving frames' would ultimately develop into the concepts of principal bundles and Cartan connections -- invariants that cover a vast amount of geometric structures and furthermore allow for explicit calculations.

This construction was further developed by T.Y.~Thomas \cite{OCG}, \cite{CT} who developed key ideas for Tractor calculus in the nineteen twenties and thirties, and S. Sasaki in 1943 \cite{Sasaki}, \cite{Sasaki2}. A major milestone was the work of N.~Tanaka \cite{OEPA} in 1979, before the seminal paper of T.N.~Bailey, M.G.~Eastwood and R.~Gover in 1994 \cite{TSB}.

Since then, there have been a series of papers by A.~\v Cap and R.~Gover \cite{TCPG}, \cite{TBIPG}, \cite{ITCCG}, \cite{ambient2}, developing a lot of the techniques that will be used in the present paper -- though those papers looked mainly at conformal geometry. Papers \cite{AHSI}, \cite{AHSII} and \cite{AHSIII}, by A.~\v Cap, J.~Slov\'ak and V.~Sou\v cek, develop similar methods in a more general setting.

Previous papers had focused on seeing the Cartan connection as a property of a principal bundle $\mathcal{P}$. But in the more recent ones, the principal bundle is replaced by an associated vector bundle, the \emph{Tractor bundle} $\mathcal{T}$, and the Cartan connection by a equivalent connection form for $\mathcal{T}$, the \emph{Tractor connection} $\overrightarrow{\nabla}$. With these tools, calculations are considerably simplified.

Here, we will start by introducing Cartan and Tractor connections for projective strusctures. We shall then relate these constructions to more standard geometric invariants -- the classes of `preferred connections', torsion-free affine connections preserving the projective structure. A few Lie algebra properties and curvature formulas will be needed to show how the Tractor connection is built up from the preferred connections.

The projective Tractor bundle $\mathcal{T}$ is of rank $n+1$, where $n$ is the dimension of the manifold. The Tractor connection preserves a volume form on $\mathcal{T}$, so we are looking at holonomy algebras contained in $\mathfrak{sl}(n+1)$. We shall first analyse the consequences of reducibility on the Tractor bundle (see Section \ref{reduc:holo}), which generate a foliation of the manifold by Ricci-flat leaves. We then look at specific cases, and show that the existence of symplectic, orthogonal, complex and hyper-complex structures on the Tractor bundle imply that the underlying manifold is projectively contact, Einstein, a $U(1)$-bundle and an $Sp(1,\mathbb{H})$-bundle respectively. Holonomies of type $\mathfrak{su}$, for instance, correspond to projectively Sasaki-Einstein manifolds. These results are summarised in table \ref{table:equi}.

\begin{table}[htbp]
\begin{center}
\begin{tabular}{|c|c|c|}
\hline
\hline
Preserved structure & Geometric structure $V$ & Equivalence? \\ 
\hline 
\hline
alternating form & Contact manifold & no \\ \hline
complex structure & $U(1)$-bundle over a complex manifold & no \\ \hline
hypercomplex structure & $Sp(1,\mathbb{H})$-bundle over a quaternionic manifold & no \\ \hline
metric & Einstein manifold & yes \\ \hline
subbundle $K \subset \mathcal{T}$ & Folliation by Ricci flat leaves & no \\ \hline
\end{tabular}
\end{center}
\caption{Tractor holonomy reduction and geometric structures}
\label{table:equi}
\end{table}

These are not equivalences, however, except in the projectively Einstein case. There are extra conditions that have to do with the rho-tensor $\mathsf{P}$, a tensor constructed bijectively from the Ricci tensor of a preferred connection. The second-order non-linear nature of $\mathsf{P}$ make these conditions somewhat subtle.

The sequel to this paper, \cite{mepro2}, will start by generating a projective cone construction, an affine, torsion-free manifold one dimension higher whose holonomy is the same as that of the Tractor connection. For this reason, we shall occasionally use the terminology for a tangent bundle connection (such as symplectic) when referring to the Tractor connection.

This cone result will allow us, using \cite{meric} and the original papers \cite{CIH} and \cite{CIH2}, to construct all examples of possible irreducible projective Tractor holonomy, and demonstrate that there are essentially no others that those described in this paper (except in the projectively Einstein case, where more variety exists). These result, to be proved in the subsequent paper, are summarised in tables \ref{table:17} and \ref{table:21}.

\begin{table}[htbp]
\begin{center}
\begin{tabular}{|c|c|c||c|c|}
\hline
\hline
algebra $\mathfrak{g}$ & representation $V$ & \ \ restrictions \ \ & algebra $\mathfrak{g}$ & representation $V$ \\ 
\hline 
\hline
& & & & \\

$\mathfrak{so}(p,q)$ & $\mathbb{R}^{(p,q)}$ & $p+q \geq 5 $ & $\widetilde{\mathfrak{g}}_2$ & $\mathbb{R}^{(4,3)}$ \\ & & & & \\

$\mathfrak{so}(n, \mathbb{C})$ & $\mathbb{C}^{n}$ & $ n \geq 5 $ & $\mathfrak{g}_2 (\mathbb{C})$ & $\mathbb{C}^7$ \\& & & & \\

$\mathfrak{su}(p,q)$ & $\mathbb{C}^{(p,q)}$ & $p+q \geq 3 $ & $\mathfrak{spin}(7)$ & $\mathbb{R}^8$ \\ & & & & \\

$\mathfrak{sp}(p,q)$ & $\mathbb{H}^{(p,q)}$ & $p+q \geq 2 $ & $\mathfrak{spin}(4,3)$ & $\mathbb{R}^{(4,4)}$ \\ & & & & \\

$\mathfrak{g}_2$ & $\mathbb{R}^7$ & & $\mathfrak{spin}(7, \mathbb{C})$ & $\mathbb{C}^8$ \\ & & & & \\
\hline
\end{tabular}
\end{center}
\caption{Projectively Einstein Holonomy algebras}
\label{table:17}
\end{table}

\begin{table}[htbp]
\begin{center}
\begin{tabular}{|c|c|c|c|c|}
\hline
\hline
algebra $\mathfrak{g}$ & representation $V$ & \ \ restrictions \ \ & manifold (local) properties \\ 
\hline 
\hline
& & & \\

$\mathfrak{sl}(n, \mathbb{R})$ & $\mathbb{R}^{n}$ & $n \geq 3 $ & Generic \\ & & & \\

$\mathfrak{sl}(n, \mathbb{C})$ & $\mathbb{C}^{n}$ & $ n \geq 3 $ & $U(1)$-bundle over a complex manifold \\& & & \\

$\mathfrak{sl}(n, \mathbb{H})$ & $\mathbb{H}^{n}$ & $ n \geq 2 $ & $Sp(1,\mathbb{H})$-bundle over a quaternionic manifold \\& & & \\

$\mathfrak{sp}(2n,\mathbb{R})$ & $\mathbb{R}^{2n}$ & $n \geq 2 $ & Contact manifold \\ & & & \\

$\mathfrak{sp}(2n,\mathbb{C})$ & $\mathbb{C}^{2n}$ & $n \geq 2 $ & Contact manifold over a complex manifold \\ & & & \\

\hline
\end{tabular}
\end{center}
\caption{Projectively non-Einstein Holonomy algebras}
\label{table:21}
\end{table}

These results are all local, avoiding issues of degneracy of the projection from the Tractor bundle to tangent bundle, and the difference between Lie algebras and their Lie groups. The manifold $M$ is always assumes to be restricted to the relevant submanifold.

The author would like to thank Dr. Nigel Hitchin, under whose supervision and inspiration this paper was crafted. This paper appears as a section of the author's Thesis \cite{methesis}. Before starting the work on Cartan and Tractor connections, we shall recall the definition of a projective structure.

\subsection{Projective structures}
A \emph{geodesic} for a manifold $M^n$ and an affine connection $\nabla$ on it is a curve $\psi : U \to M$, $U$ a subset of $\mathbb{R}$, such that
\be
\nabla_{\dot{\psi}} \dot{\psi} = 0.
\ee
An unparametrised geodesic is a curve $\psi$ such that
\be
\nabla_{\dot{\psi}} \dot{\psi} = f \dot{\psi},
\ee
for some real-valued function $f$. An unparametrised geodesic may be made into a standard geodesic by scaling $\psi$ so that $\dot{\psi}$ is replaced with $\left( \exp - {\int f d \psi } \right) \dot{\psi}$.
\begin{defi}[Projective Structure] \label{defi:pro}
A projective structure is the set of all unparametrised geodesics of a given affine connection.
\end{defi}
As we shall see, there are many affine connections preserving the same projective structure. So the projective structure is often alternately defined as:
\begin{defi}
A projective structure is an equivalence class of affine connections with the same unparameterised geodesics.
\end{defi}
For this paper, we will need to restrict attention to those affine connections that are torsion-free. This does not unduly constrain us, as
\begin{prop}
Any projective structure has a torsion-free connection compatible with it.
\end{prop}
\begin{proof}
Let $\nabla'$ be an affine connection preserving a projective structure, with torsion $\tau$. Then $\nabla = \nabla' - \frac{1}{2} \tau$ is a torsion-free connection, and if $X$ is the tangent vector of a geodesic of $\nabla'$,
\begin{eqnarray*}
\nabla_X X = \nabla'_X X - \frac{1}{2} \tau(X,X) = \nabla'_X X,
\end{eqnarray*}
so any geodesic of $\nabla'$ is a geodesic of $\nabla$.
\end{proof}
\begin{defi}[Preferred connections]
Given a manifold $M^n$ with a projective structure, a preferred connection $\nabla$ is a torsion-free affine connection preserving the projective structure.
\end{defi}

\section{Cartan and Tractor Connections}

Traditionally, since Klein, geometries were defined by a manifold $M$ and a Lie group $G$ acting transitively on $M$. The stabilizer group of any point $x \in M$ is a sub-group $P \subset G$, which changes by conjugation as $x$ varies.

From a more modern perspective, the focus has shifted to the groups $G$ and $P$, with the underlying space $M$ seen as the quotient
\begin{eqnarray*}
M = G/P.
\end{eqnarray*}
For the `flat' projective geometry, this model is $G = \mathbb{P} SL(n+1)$ and $P = GL(n) \rtimes \mathbb{R}^{n*}$. The Cartan connection is a `curved' version of these flat geometries. Given any manifold $M$, it maps the tangent space $T_M$ locally to the Lie algebra quotient,
\be
(T_M)_x \cong \mathfrak{g}/\mathfrak{p},
\ee
for all $x$ in $M$. For projective manifolds, $\mathfrak{g} = \mathfrak{sl}(n+1)$ and $\mathfrak{p} = \mathfrak{gl}(n) \rtimes \mathbb{R}^{(n*)}$. The Cartan connection solves the quivalence problem for projective structures (given the structure, there is a unique normal Cartan connection corresponding to it \cite{PGCCC} - normality is a condition similar to torsion-freeness for a Levi-Civita connection). However the Cartan connection is somewhat tricky to work with, and an equivalent construction, the Tractor connection, is often used instead. See paper \cite{TCPG} for a full study of this; but it will suffice for us to define the Tractor connection directly.

\subsection{The Tractor Connection} \label{Trac:conn}
Given a manifold $M^n$ with a projective structure, choose a preferred connection $\nabla$. Since there are no restrictions on $\nabla$ beyond the fact that it preserves the projective structure, $\nabla$ corresponds to a principal connection on $\mathcal{G}_0$, the full frame bundle of the tangent bundle $T$. Note that $\mathcal{G}_0$ has structure group $GL(n)$. One may take the contraction of the curvature $R_{hj \phantom{k} l}^{\phantom{hj} k}$ of $\nabla$ over the first and third components to get the Ricci curvature $\mathsf{Ric} \in \Gamma(T^*\otimes T^*)$. Since $GL(n)$ is reductive, $R_{hj \phantom{k} l}^{\phantom{hj} k}$ splits into a Ricci-part and a trace-free part: the Weyl curvature $W_{hj \phantom{k} l}^{\phantom{hj} k}$.

To see the relationship more clearly, we construct an equivalent tensor from $\mathsf{Ric}$, the rho-tensor $\mathsf{P}$:
\be \label{pro:rho}
\mathsf{P}_{hj} = - \frac{n}{n^2-1} \mathsf{Ric}_{hj} - \frac{1}{n^2 - 1} \mathsf{Ric}_{jh},
\ee
allowing us to write the relationship:
\be \label{pro:weyl}
R_{hj \phantom{k} l}^{\phantom{hj} k} &=& W_{hj \phantom{k} l}^{\phantom{hj} k} + \mathsf{P}_{hl} \delta_j^k + \mathsf{P}_{hj} \delta_l^k - \mathsf{P}_{jl} \delta^k_h - \mathsf{P}_{jh} \delta_l^k.
\ee

Let $L^{-n}$ be the line bundle $\wedge^n T^*$. $L^{\alpha}$ is defined to be the weight bundle $(L^{-n})^{\frac{\alpha}{-n}}$, and $T[\alpha] = T \otimes L^{\alpha}$.
\begin{defi}[Tractor Bundle]
The Tractor bundle $\mathcal{T}$ is
\be
\mathcal{T} = T[\mu] \oplus L^{\mu}
\ee
where $\mu =  \frac{n}{n+1}$.
\end{defi}
There are other `Tractor' bundles corresponding to different representations of $\mathfrak{sl}(n+1)$, the structure algebra of $\mathcal{T}$ (most notably the exterior powers of the standard representations \cite{NCKF} and the twistor representation, see \cite{TBIPG}), but we shall not need them here.

The bundle $\mathcal{A}$ of trace-free endomorphisms of $\mathcal{T}$ is consequently:
\beqa \label{A:split}
\mathcal{A} = T \oplus \mathfrak{gl}(n) \oplus T^*.
\eeqa
Here the action of $T$ natrually maps the bundle $L^{\mu}$ to $T[\mu]$ and the action of $T^*$ maps the other way. $\mathcal{A}$ is an algebra bundle; the algebraic bracket on it is given by the conditions that $[T,T] = [T^*, T^*] = 0$ and
\begin{eqnarray*}
\{ \Psi , \Pi\} &=& \Psi \Pi - \Pi \Psi, \\
\{ \Psi , X \} &=& \Psi(X), \\ 
\{ \Psi , \nu \} &=& - \Psi(\nu), \\
\{ X , \nu \} &=& X \otimes \nu + \nu (X) Id,
\end{eqnarray*}
for $\Psi, \Pi$ sections of $\mathcal{A}_0$, $X$ a section of $T$ and $\nu$ a section of $T^*$.

\begin{defi}[Tractor Connection]
The Tractor connection $\overrightarrow{\nabla}$ is given in this case by $\overrightarrow{\nabla}_X = \nabla_X + X + \mathsf{P}(X)$, or, more explicitly,
\beqa \label{tractor:formula}
\overrightarrow{\nabla}_X \left( \begin{array}{c} Y \\ a \end{array} \right) = \left( \begin{array}{c}  \nabla_X Y + Xa \\ \nabla_X a + \mathsf{P}(X,Y) \end{array} \right).
\eeqa
The dual connection on $\mathcal{T}^* = T^*[-\mu] \oplus L^{-\mu}$ is given by:
\beqa \label{dual:formula}
\overrightarrow{\nabla}_X \left( \begin{array}{c} v \\ b \end{array} \right) = \left( \begin{array}{c}  \nabla_X v - \rP(X)b \\ \nabla_X b - v\llcorner X \end{array} \right).
\eeqa
\end{defi}
The curvature of $\overrightarrow{\nabla}$ can be seen to be
\be
R^{\overrightarrow{\nabla}}_{X,Y} = \left( \begin{array}{c} 0 \\ W(X,Y) \\ CY(X,Y) \end{array} \right),
\ee
where $CY$ is the Cotton-York tensor
\be
CY_{hjk} = \nabla_h \mathsf{P}_{jk} - \nabla_j \mathsf{P}_{hk}.
\ee

\subsection{Invariance}
So far, the Tractor connection and bundle defined depend on a choice of preferred connection. How do these definitions change if we make a different choice? First, we have a triad of results about the preferred connections themselves, from \cite{TCPG} (see also \cite{methesis}):
\begin{prop} \label{line:bun}
Given a projective structure, preferred connections are in one to one correspondence with connections on any given weight bundles $L^{\alpha}$, $\alpha \neq 0$.
\end{prop}
It is easy to see that a connection on the tangent bundle must define a connection on $L^{\alpha} = (\wedge^nT^*)^{\frac{\alpha}{-n}}$ but this proposition states that the converse is also true. A consequence of this is that there are preferred connections preserving any volume form $v \in L^{-n}$; just choose the connection defined by $\nabla v = 0$. Furthermore:
\begin{prop}
The preferred connections form an affine space, modeled on $T^*$. Two preferred connections $\nabla$ and $\nabla'$ are related by a one-form $\Upsilon$ as follows:
\be
\nabla'_X Y = \nabla_X Y + \{\{X, \Upsilon\} , Y\}.
\ee
If both these connections preserve a volume form, then $\Upsilon$ is closed.
\end{prop}
And finally:
\begin{prop} \label{change:split}
The splitting of $\mathcal{T} = T[\mu] \oplus L^{\mu}$ depends on the choice of preferred connections; if we change to another preferred connection, related by $\Upsilon$, then the splitting changes by the action of $\exp{\Upsilon}$:
\beqa \label{pro:change}
\left( \begin{array}{c} Y \\ v \end{array} \right) \to \left( \begin{array}{c} Y \\ v + \Upsilon \llcorner Y \end{array} \right).
\eeqa
\end{prop}
Note that this implies that the projection $\pi^1 : \mathcal{T} \to T[\mu]$ is well-defined. This also implies (since the space of possible splittings of $\mathcal{T}$ that respect this projection is modelled on $T^*$) that the preferred connections are in one-to-one correspondence with possible splittings of $\mathcal{T}$.

These propositions are just premises to a main result of \cite{TCPG}:
\begin{theo}[Invariance]
Using the change of splitting formula of Proposition \ref{change:split} and the formula $\overrightarrow{\nabla}_X = \nabla_X + X + \mathsf{P}(X)$, the Tractor connection is defined independently of the choice of preferred connection.
\end{theo}

\section{Symplectic holonomy}
Here we shall show a strong link between symplectic holonomy and contact spaces. A reminder of the definition of a contact space:
\begin{defi}[Contact space]
A contact space is a manifold $M^{2m+1}$ with a distribution $H \subset T$ of co-dimension one that is maximally non-integrable.

Let $\theta$ be a section of the line bundle $H^{\perp} \subset T^*$. $\theta$ then defines an alternating form $d\theta$ on $H$, and a \emph{Reeb} vector field $R$ transverse to $H$ which satisfies
\be
\theta(R) &=& 1 \\
d\theta (R, - ) &=& 0.
\ee
The maximal non-integrability of $H$ is equivalent with stating that the volume form
\be
v_{\theta} = (d \theta)^m \wedge \theta
\ee
is nowhere vanishing (notice that $v_{\theta}$ is linear in the choice of $\theta$).
\end{defi}
Before discussing symplectic holonomy for $\onab$, we will need the various contact projective structures defined in \cite{fox}:
\begin{defi}[Contact projective structure]
A contact path geometry is a family of paths everywhere tangent to the contact distribution such that at each given point and each direction tangent to the contact distribution there is a unique path in the family passing through that point and tangent to that direction. A contact projective structure is a contact path geometry the paths of which are among the geodesics of an affine connection.
\end{defi}
There is an invariant of the contact projective structure, the contact torsion (essentially the torsion in the contact direction).
\begin{defi}[Contact adapted projective structure]
A contact adapted projective structure is a contact projective structure with the extra condition that the Weyl tensor $W$ of the projective structure has the property that $W_{X,Y} Z$ is a section of the contact distribution for all $X$, $Y$ and $Z$.
\end{defi}
By paper \cite{fox}:
\begin{theo}
There is an equivalence between contact projective structures with vanishing contact torsion and contact adapted projective structures.
\end{theo}
The first part of the statement purely concerns the geodesics along the contact distribution, while the second part deals with the whole projective structure. The next section will be devoted to proving that:
\begin{theo} \label{symp:theo}
If $\onab \omega = 0$ for $\omega \in \Gamma(\wedge^2 \mathcal{T})$ non-degenerate, then there is a contact adapted projective structure on the manifold, consequently a contact projective structure with vanishing contact torsion. Conversely, given any contact adapted projective structure on the manifold, there is a non-degenerate $\omega$ such that $\onab \omega = 0$.
\end{theo}
The converse comes directly from \cite{fox}, which demonstrates that the contact projective ambient construction defined therein and the ambient cone construction for projective structures (see \cite{mepro2} and \cite{oldcone}) are isomorphic - in particular, have the same holonomy, which must be that of $\onab$ as well, \cite{mepro2}. Since the contact projective ambient construction must preserve an $\omega$, so too must $\onab$. See section \ref{cartan:pro} for a discussion of these issues from the Cartan connection point of view.

\subsection{Contact adapted projective structures}
This section is devoted to proving Theorem \ref{symp:theo}. Throughout, assume that $\onab \omega = 0$ for a non-degenerate $\omega \in \Gamma(\wedge^2 \mathcal{T})$. For any bundle or element $B$, let $B^{\omega}$ the bundle that is $\omega$-orthogonal to $B$. The first step of showing that we have a contact projective structure is to identify the contact distribution:
\begin{prop}
Let $s$ be any never-zero section of $L^{\mu}$. Then the bundle $H = s^{-1}\pi^1((L^{\mu})^{\omega})$ is a contact distribution. The preferred connections are in one-to-one correspondence with connections on $L^{\omega} = T/H$. And, up to $\mathbb{R}^+$ equivalence, Reeb vector field are in one-to-one correspondence with preferred connections preserving a volume form.
\end{prop}
\begin{proof}
Notice that $H$ does not depend on the choice of $s$, as changing a bundle by a scale change $f$ send $H$ to $fH = H$. Let $\nabla$ be any preffered connection, with corresponding splitting $\mathcal{T} = T[-\mu] \oplus L^{-\mu}$. Let $t$ be any section of $L^{\mu}$ and $R^{t}$ a section of $T[\mu]$ defined by:
\be
\omega(t, R^{t}) &=& 1 \\
\omega(R^{t}, H[\mu]) &=& 0.
\ee
The second condition implies that $\omega(R^{t}, T[\mu]) = 0$, since $\omega(R^{t} , R^{t})$ is trivially zero. Differentiating the first equality gives:
\beqa
\nonumber 0 &=& \omega(\onab_X t, R^t) + \omega(t, \onab_X R^t) \\
\label{sec:cond} &=& \omega(\nabla_X t, R^t) + \omega(t, \nabla_X R^t) + \omega(t, \rP(X, R^t))\\
\nonumber &=& \nabla_X t + \omega(t, \nabla_X R^t).
\eeqa
Thus if we project $\nabla$ to being a connection on $L^{\omega} = T/H$, then if acts on $L^{\omega}$ exactly as it acts on $L^{-\mu}$. This demonstrates the one-to-one equivalence.

Let $\theta^t$ be the section of $T^*[\mu]$ defined by $\theta^t(H) = 0$ and $\theta^t(R^t) = 1$. Then the following lemma demonstrates that $H$ is a contact distribution.

If $\nabla t = 0$, then $R = t^{-1}R^{t}$ is a well defined vector field on $M$, transverse to $H$. We aim to show that it is a Reeb vector field. Since $\nabla$ preserves the volume form $t$, $\nabla$ acts on $T = t^{-1}T[\mu]$ exactly as it acts on $T[\mu]$, giving the further equalities
\be
\mathcal{T} &=& T \oplus \mathbb{R},\\
\mathcal{T} &=& T^* \oplus \mathbb{R}.
\ee
Let $\theta$ be the section of $T^*$ that is defined by $\theta(R) = 1$, $\theta(H) = 0$. To show that $H$ is a contact distribution and that $R$ is a Reeb vector field, it suffices to show:
\begin{lemm}
$d \theta$ is equal to $\omega$ on $H$ (and is hence non-degenerate) and $d\theta(R,-) = 0$.
\end{lemm}
\begin{lproof}
Let $X$ and $Y$ be sections of $H$. The conditions above imply that $\theta = \omega(t)$. Differentiating gives:
\beqa
\nonumber (\nabla_X \theta)(Y) &=& (\onab_X \theta)(Y) \\
\label{X:X} &=& \omega(\onab_X t, Y) \\
\nonumber &=& \omega(X,Y).
\eeqa
Since $\nabla$ is torsion-free, $d\theta(X,Y)$ is the skewed part of this, so $d\theta = \omega$. The same equation, after replacing $X$ with $R$ demonstrates that
\beqa \label{R:X}
(\nabla_R \theta)(Y) = 0.
\eeqa
Equation (\ref{sec:cond}) confirms that
\beqa \label{A:R}
0 = \nabla_Y t + \omega(t, \nabla_Y R) = \omega(t, \nabla_Y R) = (\nabla_Y R) (\theta)
\eeqa
Then differentiating the relation $\theta(R) = 1$ gives $(\nabla_Y \theta)(R) = 0$ and hence that $d\theta(R,Y) = 0$. Since $d\theta(R,R) = 0$ automatically,
\be
d \theta (R, -) = 0.
\ee
\end{lproof}
To show equivalence, for any given Reeb vector field, pick the preferred connection $\nabla$ such that $\nabla (R/H) = 0$. This will generate a contact one-form $\theta$, with $R$ as its Reeb vector field.
\end{proof}
So we have a well-defined contact structure on the manifold. To show that this is a contact projective structure, we need:
\begin{prop}
The geodesic defined at a point by a tangent vector $X \in \Gamma(H)$ will always remain tangent to $H$.
\end{prop}
\begin{proof}
Pick a preferred connection $\nabla$ preserving a volume form and with a Reeb vector field $R$. Let $\mu$ with tangent field $A$ be a geodesic of $\nabla$, tangent to $H$ at $p \in X$. $A = fR + X$, with $X$ a section of $H$ and $f(p) = 0$. The geodesic equation is:
\be
0= \nabla_A A &=& A(f)R + f\nabla_A R + f\nabla_R X + \nabla_X X.
\ee
Differentiating $\theta(R) =1$ and using equation (\ref{A:R}) gives $(\nabla_A R)(\theta) = 0$. Similarly, differentiating $\theta(X) = 0$ and using equations (\ref{R:X}) and (\ref{X:X}) demonstrate that $(\nabla_R X)(\theta) = (\nabla_X X)(\theta) = 0$. Thus all the terms in the geodesic equations are sections of $H$, apart from the first term $A(f)R$. That implies that $f$ must be constant along $\mu$, thus that $f=0$.
\end{proof}
This establishes that $M$ is a contact projective space. We now merely need to demonstrate the technical condition that it is adpated:
\begin{prop} \label{contact:weyl}
$W_{X,Y} Z$ is a section of $H$ for all sections $X,Y$ and $Z$ of $T$.
\end{prop}
\begin{proof}
Choose a preferred connection $\nabla$ that fixes a volume form. This gives a splitting $\mathcal{T} = T[\mu] \oplus L^{\mu} \cong T \oplus \mathbb{R}$. The Weyl curvature $W$ takes values in $\wedge^2 T^* \otimes \mathfrak{gl}(T)$. However, since $W$ is totally trace-free, $W$ actually takes values in $\wedge^2 T^* \otimes \mathfrak{sl}(T)$ - consequently $W_{X,Y} \cdot (0,1) = 0$. Since $\onab$ preserves $\omega$, $W$ must also take values in
\be
\wedge^2 T^* \otimes \mathfrak{sp}(\omega)
\ee
and thus
\be
0 &=& \omega (W_{X,Y}\cdot (Z,0), (0,1)) + \omega((Z,0), W_{X,Y}\cdot (0,1)) \\
&=& \omega (W_{X,Y}\cdot (Z,0), (0,1)),
\ee
implying that $W_{X,Y}(Z,0)$ is in $T \cap (0,1)^{\omega} = H$.
\end{proof}

\subsection{Cartan connection considerations} \label{cartan:pro}
Paper \cite{fox} deals with contact adapted projective structures, but is light on the Cartan formalism. This section is intended to recast some of those results from that point of view; it is intended for those familiar with Cartan connections in general (see \cite{PGCCC} and \cite{TCPG}, for instance).

The Lie algebras for the projective structure on a manifold are $|1|$-graded:
\be
\mathfrak{g} = \mathfrak{sl}(n+1), \ \ \mathfrak{p} = \mathfrak{gl}(n) \oplus \mathbb{R}^n
\ee
The Lie algebras for the contact projective structure are $|2|$-graded:
\be
\check{\mathfrak{g}} = \mathfrak{sp}(n+1), \ \ \check{\mathfrak{p}} = \mathfrak{csp}(n-1) \oplus (\mathbb{R}^n)^* \oplus \mathbb{R}^*,
\ee
with $\mathfrak{csp}(n-1)$ the conformal symplectic algebra $\mathbb{R} \oplus \mathfrak{sp}(n-1)$. The Lie braket on $\check{\mathfrak{p}}$ is given by the natural action of $\mathfrak{csp}(n-1)$ on $(\mathbb{R}^n)^*$ and $\mathbb{R}^*$, and by $\{\eta, \nu \} = \omega(\eta, nu)$ for a suitable choice of scaled $\omega$, when $\eta, \nu \in (\mathbb{R}^n)^*$. It can easily be seen that $\check{\mathfrak{p}} = \check{\mathfrak{g}} \cap \mathfrak{p}$ (just notice that $\check{\mathfrak{p}}$ is the subalgebra of $\check{\mathfrak{g}}$ that preserves a line).
\begin{prop}
A contact projective structure on a manifold generates a projective structure on the same manifold.
\end{prop}
\begin{proof}
The inclusion $\check{\mathfrak{p}} \subset \mathfrak{p}$ gives a bundle inclusion of the Cartan bundle $\check{\mathcal{P}} \subset \mathcal{P}$. The contact projective Cartan connection $\check{\omega}$ can then be extended to a connection $\omega$ on $\mathcal{P}$ by $P$-equivariance; this automatically transforms it into a section of $T\mathcal{P} \otimes \mathfrak{g}$, as $\mathfrak{g}$ is the span of $\check{\mathfrak{g}}$ under the action of $P$.
\end{proof}
Then what \cite{fox} effectively demonstrates is that:
\begin{prop}
$\check{\omega}$ is torsion-free if and only if $\omega$ is normal.
\end{prop}
Let us inspect what is meant by the normality of $\check{\omega}$. Recall that $\check{\mathfrak{g}}$ has a grading $\check{\mathfrak{g}}_{-2} \oplus \check{\mathfrak{g}}_{-1} \check{\mathfrak{g}}_0 \oplus \check{\mathfrak{g}}_1 \oplus \check{\mathfrak{g}}_2$ with $\check{\mathfrak{p}} = \sum_{j \geq 0} \check{\mathfrak{g}}_{j}$. The curvature tensor of $\check{\omega}$ is equivalent (using the Killing form and the identification $\check{\mathcal{P}} \times_{\check{P}} \check{\mathfrak{g}}/\check{\mathfrak{p}} = TM$) to a function $\check{\kappa} : \check{\mathcal{P}} \to \wedge^2 \check{\mathfrak{p}}^{+} \otimes \check{\mathfrak{g}}$. Normality is given by $\partial^* \kappa = 0$, where $\partial^*$ is the natural Lie Algebra co-differential. We can thus split $\kappa$ by homogeneous degree. Paper \cite{PGCCC} demonstrates that the lowest homogenity of $\kappa$ must be given by a section $\kappa_j$ taking values in $H^2(\check{\mathfrak{p}}^{+}, \check{\mathfrak{g}})$. But Konstant's version of the Bott-Borel-Weil theorem \cite{kos} gives $H^2(\check{\mathfrak{p}}^{+}, \check{\mathfrak{g}})$ as a being of homogenity two and inside $\wedge^2 \check{\mathfrak{g}}_{1} \otimes \check{\mathfrak{g}}_0$.

\begin{prop}
If $\check{\omega}$ is normal, then it is torsion-free.
\end{prop}
\begin{proof}
All other homogenity two or lower components apart from $H^2(\check{\mathfrak{p}}^{+}, \check{\mathfrak{g}})$ must vanish. But $H^2(\check{\mathfrak{p}}^{+}, \check{\mathfrak{g}}) \subset \wedge^2 \check{\mathfrak{g}}_{1} \otimes \check{\mathfrak{g}}_0$ itself is not a torsion component, and the only torsion component of higher homogeneity than two is $\wedge^2 \check{\mathfrak{g}}_{2} \otimes \check{\mathfrak{g}}_{-1}$. However this must vanish since $\check{\mathfrak{g}}_{2}$ is one-dimensional.
\end{proof}
There are various curvature conditions that allow one to establish when $\check{\omega}$ is normal, but none seem natural or informative so far. The relation between the preferred connections $\overline{\nabla}$ of the contact projective structure and $\nabla$ of the projective structure itself is worth elucidating, however. $\overline{\nabla}$ must have vanishing contact torsion; this implies that for $X$ and $Y$ sections of $H$,
\be
\overline{\nabla}_X Y - \overline{\nabla}_Y X - [X,Y] = \{X, Y \}.
\ee
Now $\{ -,- \}$ is given by a section $\nu$ of $\wedge^2 H$, depending on $\overline{\nabla}$. All the other contact torsion terms are identical to the usual torsion terms, with the exception of the $\check{T}_{1} \otimes \check{T}_{2} \otimes \check{T}_{-1}$, which is, for $R$ a section of $\check{T}_{-2}$,
\be
\overline{\nabla}_R X - \overline{\nabla}_X R - [R,X] = - \nu(\rP_{11}),
\ee
where $\rP_{11}$ is the $\check{T}_{1}^* \otimes \check{T}_{1}^*$ component of the projective contact structure. Consequently the relationship between $\overline{\nabla}$ and $\nabla$ is given by
\be
\nabla = \overline{\nabla} - \nu + \nu(\rP_{11}).
\ee

\section{Complex holonomy: covering a complex space}

Assume that $\overrightarrow{\nabla}$ preserves a complex structure $J$ on $\mathcal{T}$. Given any section $s$ of $L^{\mu}$, the vector
\be
R = s^{-1} \otimes \pi^1 J(s)
\ee
is well defined. Dividing out by the action of the one-parameter subgroup generated by $R$ gives a local projection
\be
M \to N.
\ee
\begin{theo}
The manifold $N$ has a well defined integrable complex structure $J_N$ on it.
\end{theo}
The proof of this will take several stages. The perpendicular $B^{\perp}$ to a bundle $B \subset A$ is the set of elements of $A^*$ that vanish on $B$. First notice that the perpendicular bundle $H = R^{\perp} \subset T^*$ is well defined, moreover:
\begin{lemm}
$H$ has a complex structure $J_H$ on it, derived from $J$.
\end{lemm}
\begin{lproof}
The dual tractor bundle is $\mathcal{T}^* = T[-\mu] \oplus L^{-\mu}$, where $T[-\mu]$ is the perpendicular bundle to $L^{\mu} \subset \mathcal{T}$. Then $H[-\mu] \subset T[-\mu] \subset \mathcal{T}^*$ is given by the formula
\be
H[-\mu] = (L^{\mu})^{\perp} \cap J(L^{\mu})^{\perp} = \left( L^{\mu} \oplus J(L^{\mu}) \right)^{\perp},
\ee
since $H$ is the perpendicular to $R$ and the bundle spanned by $R$ is isomorphic (up to a choice of scale) to $J(L^{\mu})$. Since the bundle $L^{\mu} \oplus J(L^{\mu})$ is automatically preserved by $J$, so is $H[-\mu]$. Call $J_H$ the restricted complex structure on $H[-\mu]$. Since $J_H$ is a section of $H[-\mu] \otimes H^*[\mu] = H \otimes H^*$, it is also a complex structure on $H$.
\end{lproof}
This $J_H$ will be the pull back of $J_N$. But before seeing this, we will need a special class of preferred connections to do calculations. These are the $R$-tangent connections.
\begin{defi}
An $R$-tangent connection $\nabla$ is a preferred connection such that in the spliting $\mathcal{T} = T[\mu] \oplus L^{\mu}$ defined by $\nabla$,
\be
J(L^{\mu}) \subset T[\mu].
\ee
\end{defi}
The name $R$-tangent comes from the fact that in this case, $J(L^{\mu})$ is isomorphic modulo a scale to the span of $R$ in $T$.
\begin{lemm}
$R$-tangent connections exist, including $R$-tangent connections that preserve a volume form.
\end{lemm}
\begin{lproof}
Being an $R$-tangent connection is a linear constraint at each point, so evidently a splitting exists with these properties. This must correspond to a preferred connection $\nabla'$. We now aim to change $\nabla'$ to an $R$-tangent connection that preserves a volume form. Let $r$ be any coordinate function with $R(r) = 1$, and let $\nu \in \Gamma(L^{-n})$ be any volume form. Then
\be
\nabla' \nu = \eta \nu,
\ee
for some one-form $\eta$. Replacing $\nu$ with $w = \exp (\int (\eta\llcorner R) dr)\nu$ gives
\be
\nabla' w = \eta' w,
\ee
where $\eta' \llcorner R = 0$. Then changing $\nabla$ by $\Upsilon = - \frac{1}{n+1} \eta'$ gives us a preferred connection $\nabla$. This connection is still $R$-tangent since the change of splitting formula of equation (\ref{change:split}) guarantees that $J(L^{\mu})$ (spanned by a scale times $R$) remains contained in $T[\mu]$. On top of this:
\be
\nabla_X w &=& \nabla'_X w + \{X, \Upsilon \} \cdot w \\
&=& -{(n+1)} (\Upsilon \llcorner X) w + \textrm{trace } (\Upsilon \llcorner X \ Id + X \otimes \Upsilon) \cdot w \\
&=& 0.
\ee
\end{lproof}

Back to the properties of $H$ and $J_H$:
\begin{lemm} \label{preserved:H}
$H$ is preserved by action of the $R$ and so is $J_H$.
\end{lemm}
\begin{lproof}
The first fact is easy to prove; it suffices to show that $\mathcal{L}_R v \in \Gamma(H)$ for any section $v$ of $H$. And this is demonstrated by:
\be
(\mathcal{L}_R v) \llcorner R &=& \big( d(v \llcorner R) + dv \llcorner R) \big)\llcorner R \\
&=& dv \llcorner R \llcorner R = 0.
\ee
To demonstrate the second fact, pick a $R$-tangent connection $\nabla$ that preserves a volume form $w$. Using $w^{\frac{\mu}{-n}}$ we may now express the splittings as
\be
\mathcal{T} &=& T \oplus \mathbb{R} \\
&\textrm{and}& \\
\mathcal{T}^* &=& T^* \oplus \mathbb{R}.
\ee
Notice that for $v$ a section of $H$, the dual Tractor connection gives, from equation (\ref{dual:formula}):
\be
\onab_R \left( \begin{array}{c} v \\ 0 \end{array} \right) &=& \left( \begin{array}{c} \nabla_R v \\ -v\llcorner R \end{array} \right) = \left( \begin{array}{c} \nabla_R v \\ 0 \end{array} \right).
\ee
Since $\onab J = 0$ and $H$ is preserved by $J$, this implies
\beqa \label{com:one}
\nabla_R Jv = J\nabla_R v.
\eeqa
Similarly, for any section $X$ of $T$:
\be
-X \llcorner J(v) &=& (\onab_X J(v)) \llcorner \left( \begin{array}{c} 0 \\ 1 \end{array} \right) \\
&=& (\onab_X v) \llcorner J \left( \begin{array}{c} 0 \\ 1 \end{array} \right) \\
&=& (\onab_X v) \llcorner \left( \begin{array}{c} R \\ 0 \end{array} \right) = (\nabla_X v) \llcorner R,
\ee
implying that
\beqa \label{com:two}
(\nabla v) \llcorner R = - J(v).
\eeqa
To prove that $\mathcal{L}_R J_H = 0$, we need to show that $\mathcal{L}_R J v = J(\mathcal{L}_R v)$. However, since $\nabla$ is torsion-free,
\be
\mathcal{L}_R J v &=& d(J(v)) \llcorner R \\
&=& (\nabla_R Jv) - (\nabla Jv) \llcorner R \\
&=& J(\nabla_R v) - v \\
&=& J(\nabla_R v + J(v)) \\
&=& J(\nabla_R v - \nabla v \llcorner R) \\
&=& J(d(v) \llcorner R) = J(\mathcal{L}_R v).
\ee
by equations (\ref{com:one}) and (\ref{com:two}).
\end{lproof}
Since $H \llcorner R = 0$, $H$ is the pull back of the bundle $TN^*$. As it has a $R$-invariant complex structure, this desends to a complex structure $J_N$ on $TN$. It suffices to show:
\begin{prop}
$J_N$ is integrable.
\end{prop}
To prove this, we need to note that:
\begin{lemm}
$J_N$ being integrable is implied by the fact that the exterior derivative $d$ on $M$ maps sections of the complex eigen-bundle $H_{\mathbb{C}}^{(1,0)}$ to sections of $H_{\mathbb{C}}^{(1,0)} \wedge T^*\mathbb{C}$.
\end{lemm}
\begin{lproof}
Since the exterior derivative commutes with pull-backs, this implies that the exterior derivative $d$ on $N$ maps sections of $(TN^*_{\mathbb{C}})^{(1,0)}$ to sections of $(TN^*_{\mathbb{C}})^{(1,0)} \wedge TN^*_{\mathbb{C}}$. Dualising this relationship implies that $TN_{\mathbb{C}}^{(1,0)}$ is closed under the Lie bracket; hence integrability.
\end{lproof}
Now we need to show that the conditions of the previous lemma do hold. Continue using the preferred connection $\nabla$ from Lemma \ref{preserved:H}. Let $s$ be the section of $\mathcal{T}^*$ corresponding to $(0,1)$ and $\tau = Js$. Since $R \subset T[\mu] \subset \mathcal{T}$, then $\tau \subset T^*[-\mu] \subset \mathcal{T}^*$. The $J$ invariance of $\onab$ gives:
\be
\onab_X Jv &=& J \onab_X v \\
\nabla_X Jv - (X \llcorner Jv)s &=& J(\nabla_X v - (X\llcorner v)s),
\ee
implying, that if $\nabla^H$ is $\nabla$ projected onto $H$ along $\tau$,
\beqa
\nabla_X v &=& \nabla_X^H v + Jv (\tau \llcorner X) \\
\nonumber &\textrm{and} & \\
\label{J:inv} \nabla^H Jv &=& J(\nabla^H v).
\eeqa
Then if $v \in \Gamma(H_{\mathbb{C}})$ we can set $v^{(1,0)} = v -iJv$. Then
\be
\nabla v^{(1,0)} &=& \nabla v - i\nabla Jv \\
&=& \nabla^H v - i\nabla^H Jv + Jv \tau + iv \tau \\
&=& \nabla^H v^{(1,0)} - i(v^{(1,0)} \tau).
\ee
The first term is a section of $T^*_{\mathbb{C}} \otimes H_{\mathbb{C}}$ and the second a section of $H_{\mathbb{C}} \otimes T^*_{\mathbb{C}}$. Since $\nabla$ is torsion-free, $d v^{(1,0)}$ is the skew-symmetrisation of this, hence a section of $H_{\mathbb{C}} \wedge T^*_{\mathbb{C}}$. This demonstrates the integrability of $J_N$, and makes $M$ locally into an $U(1)$-bundle over $N$. But is more than that; in fact, $\onab$ implies a \emph{complex projective structure} on $N$

\subsection{Complex projective structures}
A connection $\nabla$ is $R$-invariant if
\be
[R, \nabla_X Y ] = \nabla_{[R,X]} Y + \nabla_X [R,Y].
\ee
In pacticular, if $X$ and $Y$ commute with $R$, then so does $\nabla_X Y$. If a connection $\nabla$ is $R$-invariant, then so are all its curvature terms, including $\rP$, as well as any structure that $\nabla$ would preserve. Thus we can say:
\begin{defi}
A Tractor connection $\onab$ is said to be $R$-invariant if it has a preferred connection $\nabla$ that is $R$-invariant.
\end{defi}
Obviously this $\nabla$ is non-unique; it is related to all other $R$-invariant preferred connections by the action of any $R$-invariant $\Upsilon$.
\begin{defi}[$J$-preferred connections]
$J$-preferred connections are preferred connections that are $R$-invariant and $R$-tangent.
\end{defi}
If an $R$-invariant connection exists, it is easy to make if $R$-tangent as well, by first using an $\Upsilon$ doing so on any section of the projection $M \to N$, and then extending that $\Upsilon$ into $M$ by requiring it to be $R$-invariant. Thus $J$-preferred connections exist if and only if $\onab$ is $R$-invariant.

$\onab$ need not be $R$-invariant. However, we have a powerfull result if it is:
\begin{theo}
If $\onab$ is $R$-invariant, the set of $J$-preferred connections determine a complex projective structure on $N$.
\end{theo}
To prove this, we evidently have to define what we mean by a complex projective structure. $J$-preffered connections project to affine connections on $TN$, using the projection $TM/R = TN$. By equation (\ref{J:inv}), they commute with the complex structure $J_N$.

\begin{defi}[Generalised complex geodesics]
A generalised complex geodesic is a map $\psi: \mathbb{R} \to N$ such that
\begin{eqnarray*}
\overrightarrow{\nabla}_{\dot{\psi}} \dot{\psi} \in \Gamma (B),
\end{eqnarray*}
where $B$ is the bundle spanned by $\dot{\psi}$ and $J_N \dot{\psi}$. Since a real geodesic is a fortiori a generalised complex geodesic, these exist at all points, in every direction. However they are non-unique.
\end{defi}
\begin{defi}[Complex geodesics]
A complex geodesic on a complex manifold $(N,J_N, \nabla)$ is a map $\mu$ from a domain $D \subset \mathbb{C}$ to $N$ such that $\mu(U)$ is totally geodesic \cite{CGeo}, \cite{Lebrun}. They exist if the connection $\nabla$ is holomorphic -- paper \cite{CGeo} erroneously claims their existence in the general case.
\end{defi}
Obviously any curve inside a the image of a complex geodesic is a generalised complex geodesic. Note that a complex geodesic is a function $\mathbb{C} \to N$, whereas generalised complex geodesic are functions $\mathbb{R} \to N$.
\begin{prop}
All $J$-preferred connections have the same generalised complex geodesics, and, if and when they exist, the same complex geodesics.
\end{prop}
\begin{proof}
Let $\psi : \mathbb{R} \to N$ be a generalised complex geodesic for a $J$-preferred connection $\nabla$. Let $S = \dot{\psi}$ and let $\widehat{\psi}$ be any curve on $M$ that covers $\psi$ and $U = \dot{\widehat{\psi}}$. Then, by the properties of $\nabla$ and $\psi$,
\be
\nabla_U U &=& f_1 \widehat{S} + f_2 \widehat{J_N S} + f_3 R.
\ee
where $\widehat{S}$ and $\widehat{J_N S}$ are any lifts of the relevant vector fields. We may change $\nabla$ to another $J$-preferred connnection $\nabla'$ by using an $R$-invariant one-form $\Upsilon$ such that $\Upsilon \llcorner R = 0$. For this new connection,
\be
\nabla'_U U &=& \nabla_U U + \{ \{ U, \Upsilon \} U \} \\
&=& f_1 \widehat{S} + f_2 \widehat{J_N S} + f_3 R + 2 (\Upsilon \llcorner U)U.
\ee
And since $U$ is also a lift of $S$, projecting down to $N$ shows that $\nabla'_S S$ is in the span of $S$ and $J_N S$, i.e.~that $\psi$ is a generalised complex geodesic for $\nabla'$ on $N$.

Let $D \subset \mathbb{C}$ and assume that $\mu : D \to N$ is a generalised complex geodesic for $\nabla$. Then any curve $\psi$ in the image of $\mu$ is a generalised complex geodesic for $\nabla$, hence for $\nabla'$. Since $\nabla'$ commutes with $J_N$, this implies that
\be
\nabla'_{\dot{\psi}} \dot{\psi} \ \ \mathrm{and} \ \ \nabla'_{\dot{\psi}} J_N \dot{\psi}
\ee
are sections of $\mu_*(TD)$, in other words that $\mu(D)$ is a totally geodesic subspace of $N$ for $\nabla'$. Hence it is a complex geodesic for $\nabla'$, which thus shares the same complex geodesics as $\nabla$.
\end{proof}
Note that this proof did no need the $R$-invariance of $\nabla$ or $\nabla'$. Hence there may be further analogues of generalised complex geodesics even when $\onab$ is not $J$-invariant.

\begin{defi}[Complex projective structure]
A complex projective structure on $N$ is given by the complex structure $J_N$ and the generalised complex geodesics. The $J$-preferred connections are the preferred connections for this complex projective structure.
\end{defi}
By an analogous argument to that given for the real case (see Section \ref{Trac:conn}), if $\nabla$ and $\nabla'$ are two $J$-preferred connections on $N$
\begin{eqnarray*}
\widetilde{\nabla}_X Y &=& \widetilde{\nabla}'_X Y + \Upsilon^{\mathbb{C}}(X)Y + \Upsilon^{\mathbb{C}}(Y)X,
\end{eqnarray*}
with $\Upsilon^{\mathbb{C}}$ a section of $TN^*_{\mathbb{C}}$. Similarly to Proposition \ref{line:bun}, the preferred connection $\nabla$ is bijectively determined by its effect on powers of the holomorphic weight bundle
\begin{eqnarray*}
L_{\mathbb{C}}^{-n} \cong \wedge^{(n,0)} T^*_{\mathbb{C}}.
\end{eqnarray*}
We can also proceed as in the real case. Let $\widetilde{\nabla}$ be any $J$-preffered connection on $N$. The its Ricci tensor splits:
\begin{eqnarray*}
\widetilde{\mathsf{Ric}} = l_s + l_a + h_s + h_a,
\end{eqnarray*}
where
\begin{eqnarray*}
l_s(X,Y) &=& \frac{1}{4} \left( \widetilde{\mathsf{Ric}}(X,Y) + \widetilde{\mathsf{Ric}}(Y,X) - \widetilde{\mathsf{Ric}}(J_N X,J_N Y) - \widetilde{\mathsf{Ric}}(J_N Y, J_N X) \right), \\
l_a(X,Y) &=& \frac{1}{4} \left( \widetilde{\mathsf{Ric}}(X,Y) - \widetilde{\mathsf{Ric}}(Y,X) - \widetilde{\mathsf{Ric}}(J_N X,J_N Y) + \widetilde{\mathsf{Ric}}(J_N Y, J_N X) \right), \\
h_s(X,Y) &=& \frac{1}{4} \left( \widetilde{\mathsf{Ric}}(X,Y) + \widetilde{\mathsf{Ric}}(Y,X) + \widetilde{\mathsf{Ric}}(J_N X,J_N Y) + \widetilde{\mathsf{Ric}}(J_N Y, J_N X) \right), \\
h_a(X,Y) &=& \frac{1}{4} \left( \widetilde{\mathsf{Ric}}(X,Y) - \widetilde{\mathsf{Ric}}(Y,X) + \widetilde{\mathsf{Ric}}(J_N X,J_N Y) - \widetilde{\mathsf{Ric}}(J_N Y, J_N X) \right).
\end{eqnarray*}
So, $l_s$ is the $J_N$-linear symmetric component of the tensor $\widetilde{\mathsf{Ric}}$, $l_a$ the $J_N$-linear anti-symmetric component, $h_s$ the $J_N$-hermitian symmetric component and $h_a$ the $J_N$-hermitian anti-symmetric component. Then define the complex projective rho-tensor $\mathsf{P}^{\mathbb{C}}$ as
\begin{eqnarray*}
\mathsf{P}^{\mathbb{C}} &=& -\frac{l_s}{2n-2} -\frac{l_a + h_s + h_a}{2n+2}.
\end{eqnarray*}

There is also a complex projective Weyl tensor, $W^{\mathbb{C}}$. In details, this is given by
\begin{eqnarray} \label{complex:weyl} \nonumber
R_{hj \phantom{k} l}^{\phantom{hj}k} &=& (W^{\mathbb{C}})_{hj \phantom{k} l}^{\phantom{hj}k} + \big( (\mathsf{P}^{\mathbb{C}})_{hl} \delta^k_j - (\mathsf{P}^{\mathbb{C}}J_N)_{hl} (J_N)^k_j - (\mathsf{P}^{\mathbb{C}})_{jl} \delta^k_h + (\mathsf{P}^{\mathbb{C}}J_N)_{jl} (J_N)^k_h \big) \\
&& + \big( (\mathsf{P}^{\mathbb{C}})_{hj} \delta_l^k - (\mathsf{P}^{\mathbb{C}})_{jh} \delta_l^k  - (\mathsf{P}^{\mathbb{C}}J_N)_{hj} (J_N)^k_l + (\mathsf{P}^{\mathbb{C}}J_N)_{jh} (J_N)^k_l \big)
\end{eqnarray}
where $(\mathsf{P}^{\mathbb{C}}J_N)_{hj} = (\mathsf{P}^{\mathbb{C}})_{hm}(J_N)^m_j $. If we take the tensor products to be complex, this expression becomes
\begin{eqnarray*}
R_{hj \phantom{k} l}^{\phantom{hj}k} &=& (W^{\mathbb{C}})_{hj \phantom{k} l}^{\phantom{hj}k} + 2 \big( (\mathsf{P}^{\mathbb{C}})_{hl} \otimes \delta^k_j - (\mathsf{P}^{\mathbb{C}})_{jl} \otimes \delta^k_h \big) \\
&& + 2 \big( (\mathsf{P}^{\mathbb{C}})_{hj} \otimes \delta_l^k - (\mathsf{P}^{\mathbb{C}})_{jh} \otimes \delta_l^k \big).
\end{eqnarray*}
The complex Cotton-York tensor is also defined,
\begin{eqnarray*}
CY^{\mathbb{C}}(X,Y;Z) &=& + (\widetilde{\nabla}_X \mathsf{P}^{\mathbb{C}})(Y,Z) - (\widetilde{\nabla}_Y \mathsf{P}^{\mathbb{C}})(X,Z) \\
&& - i(\widetilde{\nabla}_X \mathsf{P}^{\mathbb{C}})(Y,J_N Z) + i(\widetilde{\nabla}_Y \mathsf{P}^{\mathbb{C}})(X,J_N Z).
\end{eqnarray*}
We define the Tractor bundle $\mathcal{T}_N^{\mathbb{C}}$ as the projection to $N$ of the Tractor bundle $\mathcal{T}$. In terms of the local information, a choice of $J$-preferred connections determins a splitting:
\be
\mathcal{T}_N^{\mathbb{C}} &=& \left( (T^*)^{(1,0)} \otimes_{\mathbb{C}} L_{\mathbb{C}}^{\mu_N} \right) \oplus L_{\mathbb{C}}^{\mu_N},
\ee
with $\mu_N = -\frac{m}{m+1}$, where $m = (n-1)/2$ is the complex dimension of $N$. The Tractor connection similarly projects and the expression for $\onab^{\mathbb{C}}$ on $N$ is, in local terms,
\begin{eqnarray*}
\overrightarrow{\nabla}^{\mathbb{C}}_X = \widetilde{\nabla}_X + X + \mathsf{P}^{\mathbb{C}}(X),
\end{eqnarray*}
with a complex action of $X$ and $\mathsf{P}^{\mathbb{C}}$.
\begin{rem}
There is a close connection between a change of real $J$-preferred connection on $M$, $\nabla \to \nabla'$ and the corresponding change of complex preferred connections $\widetilde{\nabla} \to \widetilde{\nabla}'$ on $N$. The first two differ by a one-form $\Upsilon$ that is zero on $R$. Then $\Upsilon$ can be made $R$-invariant by a suitable choice of isomorphisms $H \cong TN$. This makes $\Upsilon$ equivalent to a one-form $\Upsilon^{\mathbb{C}}$ on $N$, which is the one-form giving the difference between $\widetilde{\nabla}$ and $\widetilde{\nabla}'$. The converse of this is true as well.
\end{rem}

See papers \cite{CGeo} and \cite{complexprojectivestruc} for more information. The twistor results of \cite{CMEE} are also related. The sequel paper \cite{mepro2} constructs projective cone structures that tie the complex projective structure on $N$ and the real projective structure on $M$ even closer together.

\subsection{Hypercomplex holonomy}
Quaternionic and hypercomplex holonomy are treated in detail in paper \cite{mepro2}; this section just sumarises the results, since the cone construction of that paper is needed to prove them.

There are three main results:
\begin{prop}
Any $\onab$ that preserves a quaternionic structure preserves a hypercomplex structure.
\end{prop}
That mean that if $J_1'$, $J_2'$ and $J_3'$ are complex structures obeying the quaternionic relations, such that $\overrightarrow{\nabla}$ preserves their span without preserving them individualy, then $\overrightarrow{\nabla}$ must actually fix a trio of complex structures $J_1$, $J_2$ and $J_3$ (also obeying the quarernionic relations). This comes from the fact that as a consequence of \cite{mepro2}, all projective Tractor holonomies are affine holonomy algebras of torsion-free Ricci-flat cones. And $\mathfrak{sl}(1, \mathbb{H}) \oplus \mathfrak{sl}(n, \mathbb{H})$ is not a possible Ricci flat holonomy algebra by \cite{meric}.

\begin{prop}
$\overrightarrow{\nabla}$ \emph{must} be invariant in all three of the directions $R_k = s^{-1} \pi^1 J_k(s)$, $s\in L^{\mu}$.
\end{prop}
Unlike the complex case, where $R$-invariance is not required, here $R_k$-invariance is guaranteed.

\begin{theo}
Dividing out by the $R_k$, one obtains a manifold $N$ three dimensions lower, and $\overrightarrow{\nabla}$ descends to a well-defined torsion-free quaternionic connection $\nabla^N$ on $N$, i.e.~a connection with holonomy algebra contained in $\mathfrak{sl}(1, \mathbb{H}) \times \mathfrak{gl}(\frac{n-3}{4}, \mathbb{H})$.
\end{theo}
This means that ``quaternionic projective structures'' are the same thing as ``integrable quaternionic structures''. The manifold $N$ has non-canonical quaternionic structures $\{J_1^N, J_2^N, J_3^N$\}; their span forms $H_N$, a $\mathfrak{sp}(1, \mathbb{H})$ subundle of $T \otimes T^*$. The bundle $H_N$ itself is canonical for the structure, however, and is preserved by $\nabla^N$. Then $M$ itself is locally a subspace of the principle bundle for $H_N$.

\section{Orthogonal holonomy: Einstein spaces}
In this section we aim to show that $\overrightarrow{\nabla}$ preserving a metric on $\mathcal{T}$ is equivalent to the existence of an Einstein, non-Ricci-flat, preferred connection $\nabla$.

Some explanations as to what we mean by an Einstein connection in this case:
\begin{defi}
$\nabla$ is Einstein if $\mathsf{Ric}^{\nabla}$ is non-degenerate and
\begin{eqnarray*}
\nabla \mathsf{Ric}^{\nabla} = 0.
\end{eqnarray*}
\end{defi}
Notice this also implies that $\nabla \det(\mathsf{Ric}^{\nabla}) = 0$, so $\nabla$ preserves a volume form. Thus $\mathsf{Ric}^{\nabla}$ is symmetric, and $\nabla$ is the Levi-Civita connection of the `metric' $\mathsf{Ric}^{\nabla}$, meaning that $\nabla$ is an Einstein connection in the standard sense, with Einstein coefficient $1$.

\begin{prop}
If $\nabla$ is an Einstein connection, then $\overrightarrow{\nabla}$ preserves a metric $h$ on $\mathcal{T}$.
\end{prop}
\begin{proof}
Let $s \in L^{\mu}$ be a section corresponding to $\nabla$. Then in the splitting defined by $\nabla$, consider the metric
\begin{eqnarray*}
h (\left( \begin{array}{c} X \\ a \end{array} \right) , \left( \begin{array}{c} Y \\ b \end{array} \right)) = s^{-2} \left( - \mathsf{P}(X,Y) + ab \right).
\end{eqnarray*}

Note that where $\mathsf{Ric}$ is of signature $(p,q)$, $h$ is of signature $(p+1,q)$. In a more general setting, if $\mathsf{Ric} = \lambda g$ for some metric $g$ of signature $(p,q)$, then $h$ is of signature $(p+1,q)$ when $\lambda > 0$ and $(q+1, p)$ when $\lambda <0$.

Remembering the formulas for the Tractor connection, and using $s$ implicitly:
\begin{eqnarray*}
Z. h(\left( \begin{array}{c} X \\ 0 \end{array} \right), \left( \begin{array}{c} Y \\ 0 \end{array} \right)) &=& - Z. \mathsf{P}(X,Y)\\ &=& - \mathsf{P}(\nabla_Z X,Y) - \mathsf{P}(X, \nabla_Z Y) \\ &=& h (\overrightarrow{\nabla}_Z \left( \begin{array}{c} X \\ 0 \end{array} \right), \left( \begin{array}{c} Y \\ 0 \end{array} \right)) + h( \left( \begin{array}{c} X \\ 0 \end{array} \right), \overrightarrow{\nabla}_Z \left( \begin{array}{c} Y \\ 0 \end{array} \right)) \\
Z. h(\left( \begin{array}{c} X \\ 0 \end{array} \right), \left( \begin{array}{c} 0 \\ a \end{array} \right)) &=& 0 \\ &=& \mathsf{P}(Z,X) a - \mathsf{P}(Z,X) a \\ &=& h( \overrightarrow{\nabla}_Z \left( \begin{array}{c} X \\ 0 \end{array} \right) , \left( \begin{array}{c} 0 \\ a \end{array} \right)) + h(\left( \begin{array}{c} X \\ 0 \end{array} \right), \overrightarrow{\nabla}_Z \left( \begin{array}{c} 0 \\ a \end{array} \right) ) \\
Z. h(\left( \begin{array}{c} 0 \\ a \end{array} \right), \left( \begin{array}{c} 0 \\ b \end{array} \right) ) &=& 
(\nabla_Z a)b + a \nabla_Z b \\ &=& h( \overrightarrow{\nabla}_Z \left( \begin{array}{c} 0 \\ a \end{array} \right) , \left( \begin{array}{c} 0 \\ b \end{array} \right)) + h(\left( \begin{array}{c} 0 \\ a \end{array} \right), \overrightarrow{\nabla}_Z \left( \begin{array}{c} 0 \\ b \end{array} \right) ),
\end{eqnarray*}
hence
\begin{eqnarray*}
\overrightarrow{\nabla} h = 0.
\end{eqnarray*}
\end{proof}

Conversely:
\begin{prop}
If $\overrightarrow{\nabla}$ preserves a metric $h$ on $\mathcal{T}$, then there exists an Einstein preferred connection $\nabla$ on an open dense submanifold of $M$.
\end{prop}
\begin{proof}
We need first to show that $L^{\mu} \subset \mathcal{T}$ cannot degenerate for $h$, at least on an open dense subset.

Assume $h(s,s) = 0$ at $x \in M$ for some nowhere zero $s \in \Gamma(L^{\mu})$. Then
\begin{eqnarray*}
X.h(s,s) &=& 2 h(\nabla_X s, s) \\
&=& 2 h(\left( \begin{array}{c} Xs \\ \nabla_X s \end{array} \right), s) \\
&=& 2 h(\left( \begin{array}{c} Xs \\ \end{array} \right), s),
\end{eqnarray*}
and since $\left( \begin{array}{c} Xs \\ 0 \end{array} \right)$ spans an $n$-dimensional subset of $\mathcal{T}$, this quantity must be non-zero for most $X$, bar a $(n-1)$-dimensional subset of $T_x$.

Now on most points of $M$, we may define a special section $s \in \Gamma(L^{\mu})$ by requiring
\begin{eqnarray*}
h(s,s) =1.
\end{eqnarray*}
and the associated preferred connection $\nabla$ with $\nabla s = 0$. Consequently
\begin{eqnarray*}
0 &=& X.h(s,s) \\
&=& 2 h(\left( \begin{array}{c} X \\ 0 \end{array} \right), s)
\end{eqnarray*}
so $L^{\mu} \perp T[\mu]$. Moreover
\begin{eqnarray*}
0 &=& X.h(Ys,s) \\
&=& h (\left( \begin{array}{c} \nabla_X Y \\ \mathsf{P}(X,Y) \end{array} \right),s ) + h( Ys, Xs ) \\
&=& \mathsf{P}(X,Y) + h( Ys, Xs ).
\end{eqnarray*}
Hence
\begin{eqnarray*}
h (\left( \begin{array}{c} X \\ a \end{array} \right) , \left( \begin{array}{c} Y \\ b \end{array} \right)) = s^{-2} \left( - \mathsf{P}(X,Y) + ab \right),
\end{eqnarray*}
as before. As well as this,
\begin{eqnarray*}
X. \mathsf{P}(Y,Z) &=& X.h(Ys,Zs) \\
&=& h (\left( \begin{array}{c} \nabla_X Y \\ \mathsf{P}(X,Y) \end{array} \right),sZ ) + h(sY , \left( \begin{array}{c} \nabla_X Z \\ \mathsf{P}(X,Z) \end{array} \right) ) \\
&=& \mathsf{P}(\nabla_X Y,Z) + \mathsf{P}( Y, \nabla_X Z ),
\end{eqnarray*}
so
\begin{eqnarray*}
\nabla_X \mathsf{P} = 0.
\end{eqnarray*}
\end{proof}

\section{Reducible holonomy: Ricci-flatness} \label{reduc:holo}

This section will provide a description of the geometric meanings of reducible Tractor holonomy. We will not, however, fully classify this case, similar to the fact that reducible holonomy is not fully classified in the affine case. In this section, by co-volume forms, we mean elements such as
\begin{eqnarray*}
X^1 \wedge X^2 \wedge \ldots \wedge X^k
\end{eqnarray*}
where $(X^j)$ is a frame for a bundle of rank $k$.

Let $\widetilde{K} \subset \mathcal{T}$ be a rank $k \leq n$ subbundle preserved by $\overrightarrow{\nabla}$.
\begin{lemm}
On an open dense subset of the manifold, $L^{\mu}$ is not a subbundle of $\widetilde{K}$.
\end{lemm}
\begin{lproof}
This fact  is a consequence of the fact that the second fundamental form of $L^{\mu}$ is maximal, since $\overrightarrow{\nabla}$ comes from a Cartan connection (this second fundamental form is also called the soldering form).

In more details, let $\pi^1 : \mathcal{T} \to \mathcal{T} / L^{\mu} = T[\mu]$ be the quotient projection. Then the second fundamental form of $L^{\mu}$,
\begin{eqnarray*}
S: L^{\mu} \longrightarrow T^* \otimes T[\mu]
\end{eqnarray*}
is defined by
\begin{eqnarray*}
S(s)(X) &=& \pi^1 \left( \overrightarrow{\nabla}_X s \right) = sX.
\end{eqnarray*}
In consequence the image of sections of $L^{\mu}$ under $\overrightarrow{\nabla}$ span all of $\mathcal{T}$. So any bundle $\widetilde{K}$ preserved by $\overrightarrow{\nabla}$ cannot contain $L^{\mu}$ on any open set.
\end{lproof}

From now on we shall assume, by restricting to open, dense subsets of $M$, that $L^{\mu} \cap \widetilde{K} = 0$. Hence the projection $\pi^1$ is injective on $\widetilde{K}$. Given any nowhere-zero section $s$ of $L^{\mu}$, define $K \subset T$ as $s^{-1} \pi^1 (\widetilde{K})$. This bundle does not depend on a choice of $s$, as changing $s$ changes the scaling but not the bundle.

\begin{theo} \label{ric:hol}
Assume $\overrightarrow{\nabla}$ preserves $\widetilde{K}$ and $K = \pi(\widetilde{K})\subset T$. Then $K$ is an integrable, totally geodesic foliation, and there are preferred connections $\nabla$ that:
\begin{enumerate}
\item preserve $K$,
\item preserve a volume form on $K$,
\item are Ricci-flat on $K$,
\item have $\mathsf{P}^{\nabla}(-,Y) = 0$ for any section $Y$ of $K$.
\end{enumerate}
These $\nabla$ preserve a (co-)volume form on all of $T$ if and only if $\overrightarrow{\nabla}$ preserves a co-volume form on $\widetilde{K}$.
\end{theo}
Most of this section will be devoted to proving this. Given any spliting of $\mathcal{T}$, $\widetilde{K} \cap T[\mu]$ is of rank $k$ or $k-1$. If we are in the latter case, choose a frame $\{(X_1, 0), \ldots, (X_{k-1}, 0), (X_{k}, \mu)\}$ of $\widetilde{K}$. Changing this splitting by the action of $\Upsilon$ such that $\Upsilon \llcorner X_{j} = 0$, $j<k$ and $\Upsilon \llcorner X_{k} = - \mu$ gives a splitting where $\widetilde{K} \subset T[\mu]$. Let $\nabla$ be the preferred connection corresponding to this splitting.

Let $X$ and $Y$ be sections of $K$, then
\begin{eqnarray*}
\left( \begin{array}{c} Ys \\ 0 \end{array} \right)
\end{eqnarray*}
is a section of $\widetilde{K}$, for any $s \in \Gamma(L^{\mu})$. Then
\begin{eqnarray} \label{K:bundle}
\overrightarrow{\nabla}_X \left( \begin{array}{c} Ys \\ 0 \end{array} \right) = \left( \begin{array}{c} (\nabla_X Y)s + Y (\nabla_X s)\\ s \mathsf{P}(X,Y) \end{array} \right).
\end{eqnarray}
Since this must also be a section of $\widetilde{K}$, one must have $\nabla_X Y$ as a section of $K$, and consequently $[X,Y] = \nabla_X Y - \nabla_Y X$ is a section of $K$. Hence
\begin{prop}
$K$ is integrable and totally geodesic.
\end{prop}
If one were choose $X$ as any section of $T$ rather than $K$ in Equation (\ref{K:bundle}), one sees that $\nabla$ preserves $K$ and
\begin{eqnarray*}
\mathsf{P}(-,Y) = 0,
\end{eqnarray*}
since $K$ has no $L^{\mu}$ component.

\begin{rem}
Note that as a consequence of this, $\mathsf{P}$ is zero on $K \otimes K$, hence $\mathsf{Ric}$ is zero on this foliation as well. Since $K$ is preserved by $\nabla$ in all directions, $\mathsf{Ric}^K = \mathsf{Ric}^M|_{K \otimes K}$ (this may be seen directly by taking a frame of $K$ and extending to a frame of $T$). In other words, the leaves of the foliations $K$ are Ricci-flat under the connection $\nabla$ restricted to these leaves.
\end{rem}
\begin{lemm}
We may choose $\nabla$ so that it preserves a co-volume form on $K$.
\end{lemm}
\begin{lproof}
Since $\nabla |_K$ is Ricci-flat, it must preserve a co-volume form $\tau$ \emph{along} $K$. Thus
\begin{eqnarray*}
\nabla \tau = \omega \otimes \tau,
\end{eqnarray*}
where $\omega$ is a one-form with $\omega(K) = 0$. Now $\{\Upsilon, X \}$ acts on $\tau$ by taking the trace of the first $k$ components; or, in other words,
\begin{eqnarray*}
\{\Upsilon, X \}. \tau &=& - \left( \Upsilon \otimes X + \Upsilon(X) Id \right) \tau \\
&=& -\Upsilon(\tau) \wedge X - k(\Upsilon(X)) \tau.
\end{eqnarray*}
In other words, if we change preferred connections from $\nabla$ to $\nabla'$ by the choice of
\begin{eqnarray*}
\Upsilon = - \frac{1}{k} \omega,
\end{eqnarray*}
then
\begin{eqnarray*}
\nabla' \tau = \omega \otimes \tau - \frac{k}{k} \omega \otimes \tau = 0.
\end{eqnarray*}
Since $\Upsilon(K) = 0$, then by Equation (\ref{pro:change}), $\nabla'$ still determines a splitting with $\widetilde{K} \subset T[\mu] \subset \mathcal{T}$.
\end{lproof}
\begin{prop}
There is a relationship between the holonomy of $\overrightarrow{\nabla}$ and the properties of this $\nabla'$: $\overrightarrow{\nabla}$ preserves a co-volume form on $\widetilde{K}$ if and only if $\nabla'$ preserves a (co-)volume form on $T$.
\end{prop}
\begin{proof}
From equation (\ref{tractor:formula}) and since $\mathsf{P}^{\nabla'}(-,Y) = 0$ for any section $Y$ of $K$, $\overrightarrow{\nabla}$ acts on $\widetilde{K}$ in the same way that $\nabla'$ acts on $K[\mu]$. If $\nabla'$ preserves a nowhere zero section $s$ of $L^{\mu}$, then $\overrightarrow{\nabla}$ preserves $s^k \tau$ on $\widetilde{K}$.

Conversely, if $\overrightarrow{\nabla}$ preserves a co-volume form $\widetilde{\tau}$ on $\widetilde{K}$, then
\begin{eqnarray*}
\widetilde{\tau} = t \tau
\end{eqnarray*}
for $t$ some nowhere-zero section of $L^{k \mu}$. Then
\begin{eqnarray*}
0 &=& \nabla' \widetilde{\tau} \\
&=& \nabla' t \tau \\
&=& (\nabla' t ) \tau + t(\nabla' \tau) \\
&=& (\nabla' t ) \tau.
\end{eqnarray*}
Hence $\nabla' t = 0$.
\end{proof}

\begin{cor}
Theorem \ref{ric:hol} clearly has a converse: let $\nabla$ be a preferred connection with a preserved totally geodesic integrable foliation $K$ such that $\mathsf{P}^{\nabla}(Y,-) = 0$ for a section $Y$ of $K$. Then $\overrightarrow{\nabla}$ preserves a subbundle $\widetilde{K}$ of $\mathcal{T}$. If there exists preferred connections with these properties which preserve co-volume forms on $K$ and a (co-)volume form on $T$, then $\overrightarrow{\nabla}$ preserves a co-volume form on $\widetilde{K}$.
\end{cor}
As a consequence of this, if $\widetilde{K}$ is a rank $n$ bundle, then $K = T$, and there exists a Ricci-flat preferred connection $\nabla$ on $M$. Since it is Ricci-flat, it must preserve a volume form, hence:
\begin{cor}
If $\overrightarrow{\nabla}$ preserves a rank $k = n$ bundle $\widetilde{K}$, it always preserves a volume form on $\widetilde{K}$.
\end{cor}

Notice that since the rho-tensor of $\nabla$ is zero on $K$, as is the rho-tensor of $\nabla|_K$, the tractor connection of $K$ is a restriction of that of $M$:
\begin{eqnarray*}
\overrightarrow{\nabla}^K_X \left( \begin{array}{c} Ys^{\nu} \\ t^{\nu} \end{array} \right) = \overrightarrow{\nabla}_X \left( \begin{array}{c} Ys \\ t \end{array} \right)
\end{eqnarray*}
whenever $X$ and $Y$ are sections of $K$, and $\nu = \frac{\mu_K}{\mu_M} = \frac{n(k+1)}{(n+1)k}$.

There is another useful characterisation in the `nearly irreducible' case, where $n = k$:
\begin{theo} \label{prohol:riccif}
If $\overrightarrow{\nabla}$ preserves a bundle $K$ of rank $n$ and acts irreducibly on $K$ then the holonomy algebra of $\overrightarrow{\nabla}$ is
\begin{eqnarray*}
\overrightarrow{\mathfrak{hol}} = \mathfrak{hol}^{\nabla} \oplus T \ \ \textrm{or} \ \ \overrightarrow{\mathfrak{hol}} = \mathfrak{hol}^{\nabla},
\end{eqnarray*}
where $\mathfrak{hol}^{\nabla}$ is the affine holonomy algebra of the Ricci-flat preferred connection $\nabla$ on $M$. The Lie bracket is given by the standard one on $\mathfrak{hol}^{\nabla}$, the trivial one on $T$, and action of $\mathfrak{hol}^{\nabla}$ on $T$ in cross terms.
\end{theo}
\begin{proof}
Remember the algebra bundle splitting from equation (\ref{A:split}):
\begin{eqnarray*}
\mathcal{A} = T^* \oplus \mathfrak{gl}(T) \oplus T.
\end{eqnarray*}
In the splitting given by $\nabla$, $T[\mu] = \widetilde{K}$ is preserved by $\overrightarrow{\nabla}$, thus there can be no $T^*$ component to the holonomy of $\overrightarrow{\nabla}$. As $\overrightarrow{\nabla}$ and $\nabla$ act identically on $T[\mu]$, the $T \otimes T^*$ component of the holonomy of $\overrightarrow{\nabla}$ must be the affine holonomy of $\nabla$. Then given the conditions on $\overrightarrow{\nabla}$, $\mathfrak{hol}^{\nabla}$ must act irreducibly on $T[\mu]$.

Then the algebra $\mathfrak{hol}^{\nabla} \oplus T$ decomposes into two pieces, $\mathfrak{hol}^{\nabla}$ and $T$, under the action of $\mathfrak{hol}^{\nabla}$. In other words, if the holonomy of $\overrightarrow{\nabla}$ has any $T$ component, it has the full $T$.

In actual fact, (see \cite{methesis}) $\overrightarrow{\mathfrak{hol}} = \mathfrak{hol}^{\nabla}$ if and only if $M$ is a projective cone in the sense of paper \cite{mepro2}.
\end{proof}

There is no complementary foliation to $K$ and the condition $\mathsf{P}(-,Y) = 0$ is a second order non-linear differential one; consequently it is hard to understand exactly what restrictions they impose on the projective structure. A pair of examples from the author's thesis \cite{methesis} suffice to show that these restrictions are geometrically not that strong, even when the various dimensions or co-dimensions are low.
\begin{prop}
The condition $\mathsf{P}(-,Y) = 0$ is truly a restriction on the rho-tensor; one may have connections $\nabla$ with this property where $\mathsf{Ric}^\nabla (-,Y) \neq 0$, even when $K$ is of co-dimension one.
\end{prop}
and
\begin{prop}
Assume $k \geq 3$. Then $\nabla$ restricted to one leaf of $K$ may be flat, even if it is non-flat and with maximal holonomy $\mathfrak{sl}(k)$ when restricted to a different leaf of $K$. This result remains valid if $\nabla$ preserves a volume form or not, is Ricci-flat or not, and whatever the codimension of $K$ is.
\end{prop}
These results can then be generalised to a wide variety of varying holonomy groups. So it seems that the condition $\mathsf{P}(-,Y) = 0$ is not enough to pin down the geometry in any significant way.

\end{document}